\documentclass[letterpaper, 10 pt, conference]{ieeeconf}  
\IEEEoverridecommandlockouts    
\usepackage{cite}
\newcommand{\eat}[1]{}
\usepackage{booktabs}
\usepackage{color}
\usepackage [autostyle, english = american]{csquotes}
\MakeOuterQuote{"}
\usepackage{xcolor}

\usepackage{multicol}

\usepackage{mathtools}
    
\usepackage{amsmath}
\usepackage{amstext}
\usepackage{amssymb}
\usepackage{amsfonts}
\usepackage{float}

\usepackage{amsthm}  

\usepackage{subfig}
\usepackage{caption}
\usepackage{todonotes}

\makeatletter

\newcommand{\Rmnum}[1]{\expandafter\@slowromancap\romannumeral #1@}
\usepackage{savesym}

\usepackage{algorithm}
\usepackage{algorithmicx} 
\usepackage{algpseudocode} %
\savesymbol{AND}
\usepackage[group-separator={,},group-minimum-digits={3}]{siunitx}

\usepackage{graphicx} 
\usepackage{epsfig} 

\usepackage{times} 
\usepackage{amsmath} 
\usepackage{amssymb}  

\makeatletter
\let\NAT@parse\undefined
\makeatother
\usepackage{hyperref}
\hypersetup{
   colorlinks=true,
    linkcolor= blue,
    allcolors=blue,
    citecolor = blue,
    filecolor=black,      
    urlcolor=blue,
    }
\usepackage{mathrsfs}
\title{\LARGE \bf
A Hysteretic Q-learning Coordination Framework for Emerging Mobility Systems in Smart Cities

}
\eat{A Learning-Based Approach for Optimal Coordination of Connected and Automated Vehicles}
\author{Behdad Chalaki, \emph{IEEE Student Member}, Andreas A. Malikopoulos, \emph{IEEE Senior Member}%
\thanks{This research was supported in part by ARPAE’s NEXTCAR program under the award number DE- AR0000796 and by the Delaware Energy Institute (DEI). This support is gratefully acknowledged.}%
\thanks{The authors are with the Department of Mechanical Engineering, University of Delaware, Newark, DE 19716 USA (emails: \texttt{\{bchalaki;andreas\}@udel.edu}).}}

\begin{document}

\maketitle
\thispagestyle{empty}
\pagestyle{empty}

\begin{abstract}
Connected and automated vehicles (CAVs) can alleviate traffic congestion, air pollution, and improve safety. In this paper, we provide a decentralized coordination framework for CAVs at a signal-free intersection to minimize travel time and improve fuel efficiency. We employ a simple yet powerful reinforcement learning approach, an off-policy temporal difference learning called Q-learning, enhanced with a coordination mechanism to address this problem. Then, we integrate a first-in-first-out queuing policy to improve the performance of our system. We demonstrate the efficacy of our proposed approach through simulation and comparison with the classical optimal control method based on Pontryagin's minimum principle.  

\end{abstract}

\section{Introduction}
\subsection{Motivation}
\PARstart{O}{ver} the last decade, the growing population in urban areas without increasing the road capacities has led to traffic congestion, increasing delays, and environmental concerns \cite{Schrank2019}. However, congestion might become more severe in traffic scenarios, such as intersections. When it comes to traffic safety, intersections are one of the main contributors to traffic injuries among different traffic scenarios; an average of one-quarter of traffic fatalities and roughly half of all traffic injuries are attributed to intersections \cite{FHWA1}.

The introduction of communication technologies along with computational capabilities into \textit{connected and automated vehicles} (CAVs) has the potential to revolutionize the transportation system. Through these advancements, our transportation system transitions to an \textit{emerging mobility system}, in which CAVs can make better operational decisions leading to significant reductions of energy consumption, greenhouse gas emissions, travel delays, and improvements to passengers safety \cite{Klein2016a,Melo2017a,ersal2020connected,Wadud2016}. 

\subsection{Related Work} 
After the seminal work of Levine and Athans \cite{Levine1966,Athans1969} on safely coordinating vehicles at merging roadways, several research efforts have explored the benefits of coordinating CAVs in traffic scenarios, such as urban intersections \cite{Malikopoulos2017,hult2018optimal,bichiou2018developing,Lee2012a,chalaki2020TCST}, merging roadways\cite{Ntousakis:2016aa,xiao2019decentralized,Rios-Torres2017b}, and speed reduction zones \cite{Malikopoulos2018c} to eliminate congestion in a transportation network while preserving safety by using classical control techniques.  A compendious survey of the research efforts reported in the literature to date in control and coordination of CAVs using classical control approaches is provided in \cite{guanetti2018control} and \cite{Rios-Torres2017}.

Although classical control techniques are effective for some traffic control tasks, the complexity and need for a complete knowledge of the system dynamics make autonomous driving a notoriously difficult problem to address in the real systems. On the other hand, the evolution of processing power and generation of a massive amount of data have paved the way for reinforcement learning (RL) techniques to emerge as an alternative method for traffic control. RL approaches are used when an agent learns from interaction with an environment without requiring the complete models of environment. These approaches employ the Markov decision process framework to define the interaction between a learning agent and its environment \cite{sutton2018reinforcement}. Kiumarsi \textit{et al.} \cite{kiumarsi2017optimal} surveyed various RL-based techniques to solve optimal control problems in real-time using data measurement along the system trajectories.

An off-policy temporal difference learning called Q-learning Q is one of the simplest and most promising RL methods introduced by Watkins \cite{watkins1989learning} in $1989$. Since then, numerous studies have been reported in the literature to employ Q-learning in various transportation applications. El-Tantawy and Baher \cite{el2010agent} presented a Q-learning-based traffic signal control with different state representations and tested the effectiveness of their model on a real-world multi-phase intersection. A group of papers employed a Q-learning-based approaches in ramp-metering control \cite{ji2009optimal,jacob2005integrated,davarynejad2011motorway,ivanjko2015ramp}. Wang \textit{et al.} \cite{wang2019q} used a Q-learning approach for smart lane-changing maneuvers in a mixed-traffic scenario at freeways to improve the travel time and traffic flow. They showed that their proposed approach outperformed the baseline scenario for market penetration rate below 60\%. Ngai and Yung \cite{ngai2011multiple} established a multiple-goal framework for overtaking CAVs based on Q-learning and double-action Q-learning. They decomposed the overtaking problem into several different goals that the learning agent needs to achieve. They used a double-action Q-learning formulation for the collision avoidance goal and employed Q-learning for other goals. Wu \textit{et al.} \cite{wu2019dcl} presented a decentralized learning-based framework for autonomous intersection management. The authors modeled CAVs crossing the intersection as a \textit{multi-agent Markov decision process} (MAMDP), in which CAVs cooperate to minimize the intersection delay, and solved it through Q-learning. To mitigate the "curse of dimensionality" and environment nonstationarity, they decomposed the state space of the system for each agent into independent and coordinated parts. The authors updated the corresponding Q-values for those parts separately for each CAV.

\eat{However, in large problems with many state-action pairs, to avoid Bellman's "curse of dimensionality," one can approximate the Q-function with a deep-neural-network, commonly known as deep-Q-network. This method, which was first introduced by Mnih \textit{et al.} \cite{mnih2013playing} to learn control policies for computer games using only raw pixels, has been explored by researchers in various transportation applications. 
Seliman \textit{et al.} \cite{seliman2020automated} proposed a control strategy based on a Deep-Q-Network (DQN) Reinforcement learning (RL) for a single CAV to navigate through a congested lane-reduction zone. Isele and Cosgun \cite{isele2017go} established a DQN framework for navigating a single CAV through a signal-free intersection and compared its performance with the baseline policy in a traffic simulator. Tram \textit{et al.} \cite{tram2018learning} added a recurrent layer to the DQN for navigating a single CAV among human-driven vehicles with unknown intentions at an intersection. Instead of finding a policy with continuous control input, the authors used short-term goals as an action of the agent, in which the agent learns to keep a set speed, yield for a crossing car, or stop before the intersection. Several other efforts have considered using DQN in power-train optimization of hybrid electric vehicles \cite{he2019deep}, eco-routing \cite{ma2019meta}, and  overtaking maneuver \cite{ronecker2019deep}. A survey of the state of the art on deep learning applications on control of CAVs is given in \cite{kuutti2020survey}.}

\subsection{Contribution of This Paper}
Although there have been several research efforts reporting on Q-learning-based frameworks for different transportation applications, to the best of our knowledge, no paper has reported work on a decentralized RL-based coordination framework for CAVs at an intersection intending to minimize energy consumption and improve traffic throughput. 
In this paper, we establish a decentralized coordination framework for CAVs at a signal-free intersection to minimize travel time and improve fuel efficiency. We formulate the problem by employing a well-known RL approach enhanced with a coordination mechanism called a hysteretic Q-learning, in which two learning rates are considered. Additionally, we integrate a first-in-first-out (FIFO) queuing policy in our hysteretic Q-learning framework to improve the performance of our system. We show our proposed approach's effectiveness through simulation and comparison with the classical optimal control method based on Pontryagin's minimum principle. The contributions of this paper are: (1) the development of a hysteretic Q-learning optimal framework to coordinate CAVs at a signal-free intersection aimed at decreasing both travel time and fuel consumption of each CAV; (2) integrating FIFO queuing policy into our Hysteretic Q-learning optimal framework; and (3) comparison of the proposed framework with the benchmark solution from the classical control method based on Pontryagin's minimum principle.

 \subsection{Comparison with Related Work}
 The proposed framework advances the state of the art in the following ways. First, rather than considering a single agent in the RL framework \cite{isele2017go,seliman2020automated,tram2018learning}, we propose a decentralized multi-agent framework with $100\%$ penetration rate of CAVs.
Second, in contrast to \cite{wu2019dcl,isele2017go,tram2018learning}, we incorporate energy consumption minimization into our framework in addition to traffic throughput improvement while ensuring both lateral and rear-end safety through our combined hysteretic Q-learning with FIFO framework. Third, in contrast to the research efforts reported in the literature to date, we compare the results of our proposed framework with the classical control method based on Pontryagin's minimum principle. 

\subsection{Organization of This Paper}
The rest of the paper is structured as follows. 
In Section \Rmnum{2}, we introduce the modeling framework. After introducing some preliminary materials, we provide our hysteretic Q-learning framework. Then, we present our combined hysteretic Q-learning with FIFO framework. We present the simulation framework in Section \Rmnum{3} and the corresponding results in Section \Rmnum{4}. Finally, we present concluding remarks and some discussion for a future research direction in Section \Rmnum{5}.

\section{Problem Formulation} 
\begin{figure}
\centering
\includegraphics[width=0.99\linewidth]{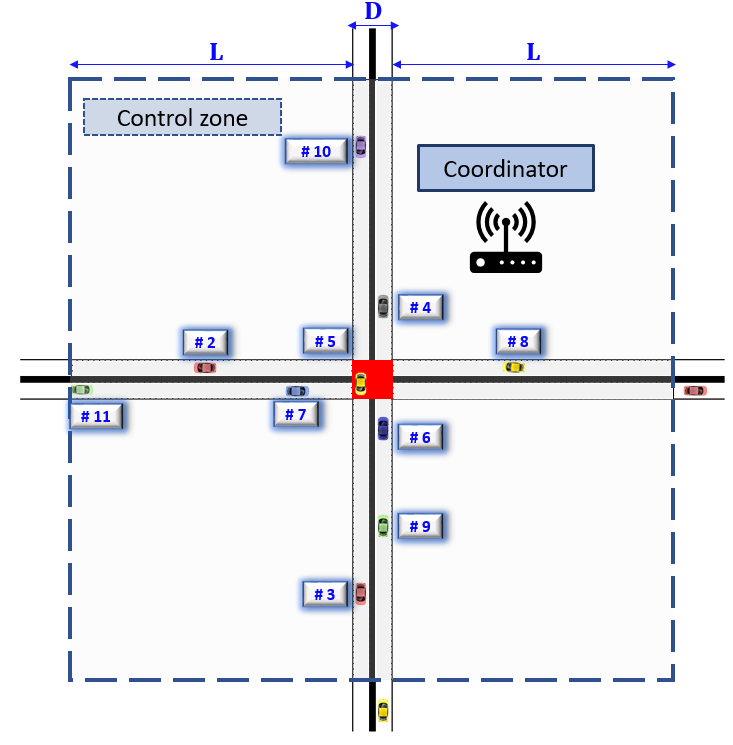}
\caption{An intersection with a coordinator communicating with CAVs inside the control zone.}
\label{fig:1}%
\end{figure}

We consider a signal-free intersection (see Fig. \ref{fig:1}), which includes a \textit{coordinator} that stores information about the intersection's geometric parameters and CAVs' information. The coordinator does not make any decision, and it only acts as a \textit{\mbox{{database}}} among the CAVs. The intersection includes a \textit{{contol zone}} in which the coordinator can communicate with the CAVs. We assume, there are no errors or delays in the vehicle-to-vehicle and vehicle-to-infrastructure communication. Although this is a strong assumption, it is relatively straightforward to relax this assumption as long as the noise or delays are bounded. We call the area inside the control zone where lateral collisions may occur \textit{{merging zone}}. The distance from the entry of the control zone to the entry of the merging zone is $L\in\mathbb{R}_{>0}$, and it is assumed to be the same for all entry points. The length of the merging zone is denoted by $D\in\mathbb{R}_{>0}$ (Fig. \ref{fig:1}). We limit our analysis to the cases where left/right turns and lane-changing maneuvers are not allowed.

The problem is formulated as a multi-agent Markov decision process (MAMDP) $\left\langle n,~\mathcal{S},~\mathcal{U},~P,~R,~\gamma\right\rangle$, where $n\in\mathbb{N}$ is total number of CAVs, $\mathcal{S}\coloneqq  \times^{n}_{i=1} \mathcal{S}^{i}$ is a finite set of states of all CAVs, $\mathcal{U}\coloneqq  \times^{n}_{i=1} \mathcal{U}^{i}$ is the joint action space, where $\mathcal{U}^{i},~ i\in\{1,2,\ldots,n\}$ is the finite set of actions of CAV $i$, $P\coloneqq \mathcal{S}\times\mathcal{U}\times \mathcal{S}\rightarrow[0,1]$ is the state transition probability which defines the transition probability between states, $R\coloneqq \mathcal{S}\times\mathcal{U}\times \mathcal{S}\rightarrow\mathbb{R}$ is the reward function for all CAVs, $R^i$ is the reward function for CAV $i$, and $\gamma \in [0,1]$ is a discount factor. 

Next, we briefly explain different approaches to formulate the Q-learning updates along with advantages and disadvantages of each approach. For CAV $i$, we use $s^i_k$ and $u^i_k$ to denote the state and action that CAV $i$ takes at time step $k\in\mathbb{N}$, respectively. Taking action $u^i_k\in\mathcal{U}^i$, CAV $i$ transitions from $s^i_k\in\mathcal{S}^i$ to $s^i_{k+1}\in\mathcal{S}^i$ and receives the reward $r^i_k=R^i(s^i_k,u^i_k)$. In the centralized update rule, the multi-agent system is viewed as a whole and is solved as a single-agent learning task, in which there is only a single Q-function. The update rule is
\begin{align}\label{eq:centralized:update}
&Q(s_k,u^1_k,\dots,u^n_k) \leftarrow (1-\alpha)Q(s_k,u^1_k,\dots,u^n_k)\nonumber\\ +&\alpha \left[r_k + \gamma \max\limits_{u_{k+1}^1,\cdots,u_{k+1}^n} Q(s_{k+1},u^1_{k+1},\dots,u^n_{k+1})\right],
\end{align}
where $s_k\in\mathcal{S}$ is the state of the system (collection of states of all CAVs), $r_k=R(s_k,u^1_k,\ldots,u^n_k)$ is the total cost incurred on the system at time step $k\in\mathbb{N}$, and $\alpha \in (0,1]$ is the learning rate. Although, theoretically, this approach converges with probability $1$ to the optimal action-value function, it does not scale well when the number of agents is increasing as the size of Q-table grows exponentially.
 
In the decentralized framework, each CAV is an \textit{independent learner} (IL) with a corresponding Q-function. 
The update rule for CAV $i$ is
\begin{align}\label{eq:decentralized:update}
&Q^i(s_k^i,u_k^i) \leftarrow Q^i(s_k^i,u_k^i)\nonumber\\& +\alpha \left[r^i_k + \gamma \max\limits_{u_{k+1}^i} Q^i(s_{k+1}^i,u_{k+1}^i) -Q^i(s_k^i,u_k^i)\right],
\end{align}
where $s_{k+1}^i\in\mathcal{S}^i$ is the state of CAV $i$ at time step $k+1$.
This approach has a smaller size Q-table than the centralized approach, and by increasing the number of agents, the Q-table's size does not grow exponentially. However, one of the drawbacks of this method is the lack of any coordination mechanism.

In our problem, CAVs need to coordinate to cross the intersection safely. Without a coordination mechanism, a CAV may select an optimal action, but it gets penalized due to the sub-optimal actions of other CAVs. It has been shown in \cite{wu2019dcl,meng2020adaptive,matignon2007hysteretic}, that using decentralized learning in a multi-agent framework with interacting agents leads to more oscillation in the learned policy and poorer performance compared to the centralized approach. In addition, since all CAVs are learning synchronously, the environment is not stationary anymore from the perspective of any CAV. Since past actions of some CAVs may affect the current behavior of other CAVs, the system is not Markovian. The latter implies that convergence is not guaranteed for every single CAV \cite{matignon2007hysteretic}.

\subsection{Hysteretic Q-learning Framework}
 Matignon \textit{et al.} \cite{matignon2007hysteretic} first presented the hysteretic Q-learning approach to incorporate coordination among ILs in a decentralized RL framework by including two learning rates. They showed that by incorporating two learning parameters $\alpha$ and $\beta$, without affecting the Q-table size, the coordination among IL agents could be achieved. In addition, the performance of the system is as good as the centralized approach of multi-agent RL. The update rule for the hysteretic Q-learning is 
 
 \begin{align}\label{eq:TDError}
&\delta \leftarrow r_k^i + \gamma \max\limits_{u_{k+1}^i} Q^i(s_{k+1}^i,u_{k+1}^i) -Q^i(s_k^i,u_k^i),
\end{align}
where
 \begin{equation}\label{eq:HystUpdate}
    Q^i(s_k^i,u_k^i) = \begin{cases}
             Q^i(s_k^i,u_k^i) + \alpha \delta,& \text{if}~ \delta\geq0,
            \\
            Q^i(s_k^i,u_k^i) + \beta \delta, &\text{otherwise,}
          \end{cases}
\end{equation} 
where $\delta$ is a temporal difference (TD) error and $\beta<\alpha \in (0,1]$. By using smaller learning rate when TD error is negative, the update results in a slower degradation of Q-value (hysteresis) associated with positive past experience. For instance, due to the sub-optimal actions of other CAVs in the environment, CAV $i$ may get penalized by doing action $u^i$ at state $s^i$, for which it received a reward in the past. In this case, the effects of this penalty on Q-value of agent $i$ should be less important.

\subsection{Main Elements of Proposed Framework}
In the problem we consider in this paper, we adopt the hysteretic Q-learning formulation to update the Q-tables. To the best of our knowledge, this is the first analysis of a transportation system using hysteretic Q-learning. In our framework, each CAV is an IL agent with a unique assigned index. Next, we present the main elements of our approach, including states, actions, and rewards.
\subsubsection{State space}
At time step $k$, we consider that CAV $i$ partially observes the system, and its state is $s^i_k\coloneqq\left\langle p^i_k,~v^i_k,~\mathcal{X}^{\text{i,rear}}_k,~\mathcal{P}^{\text{i,lat}}_k \right\rangle$, where $p^i_k$ and $v^i_k$ are its position and speed, respectively; $\mathcal{X}^{\text{i,rear}}_k\coloneqq\left\langle p^j_k,~v^j_k\right\rangle$ consists of the position and speed of CAV $j$, which is immediately ahead of CAV $i$ (e.g. CAV \#$6$ is immediately ahead of CAV \#$9$ in Fig. \ref{fig:1}); and $\mathcal{P}^{\text{i,lat}}_k$ consists of the position of the three closest vehicles to the exit of merging zone that have a potential of lateral collision with CAV $i$ (e.g. CAV \#$7$, \#$8$, and \#$10$ are three closest vehicles to the exit of merging zone for CAV \#$9$ in Fig. \ref{fig:1}). \eat{This state space was designed with the real-world implementation in mind, and could conceivably be implemented through vehicle-to-vehicle (V2V) communication between the CAVs, or vehicle-to-infrastructure (V2I) communication between CAVs and the coordinator.}

\subsubsection{Action space}
CAV $i$ has to choose action $u^i_k$ at time step $k$ which is acceleration/deceleration from a discrete bounded set $\mathcal{U}^i$ with lower bound $u^{i,\min}$ and  upper bound $u^{i,\max}$ which correspond to the minimum and maximum allowable control input of CAV $i$, respectively. Without loss of generality, we do not consider variation among CAVs' maximum and minimum control input. To this end, we set $u^{i,\min}=u^{\min}$ and $u^{i,\max}=u^{\max}$. In order to choose all actions in all states with nonzero probability and balance between exploration and exploitation, we employ the epsilon-greedy algorithm with a linear decay as follows

 \begin{align}\label{eq:decay}
 \rho&= \max\left\{\dfrac{\text{total episodes}-\text{current episode}}{\text{total episodes}},0\right\},\\
            \epsilon &=(\epsilon_{i}-\epsilon_{f})\rho + \epsilon_{f},\label{epsilon}
\end{align}
where $\rho$, $\epsilon_{i}$, and $\epsilon_{f}$ are decay rate, initial and final ratio of exploration, respectively.
In a RL framework, each episode represent a simulation, in which there is a corresponding epsilon found from \eqref{epsilon}. The corresponding epsilon determines the probability that an agent takes a random action at each episode. The epsilon found from \eqref{epsilon} is bounded between initial and final ratio of exploration.
 \begin{equation}\label{eq:eps-greedy}
    u_k^i = \begin{cases}
             \arg\max\limits_{u_k^i} Q^i(s_k^i,u_k^i),& \text{with probability}~ 1-\epsilon,
            \\
            \text{random action},& \text{with probability}~ \epsilon,
          \end{cases}
\end{equation} 
where $\epsilon$ is a small positive number. Employing epsilon-greedy with a linear decay results in more exploration at the earlier episodes and less at the final episodes which can improve the performance of the framework.
It is worth mentioning that by choosing a same value for $\epsilon_{i}$ and $\epsilon_{f}$, the algorithm simplifies to the epsilon-greedy algorithm without decay.

\subsubsection{Rewards}
CAV $i$ takes an action $u^i_k$ at time step $k$, transitions from state $s^i_k$ to the new state $s_{k+1}^i$, and receives a reward (or penalty) $r_{k}^i$ based on the following multi-objective cost:
 \begin{equation}\label{eq:rewardtotal}
r_k^i=w_1\cdot r^i_{\text{fuel}}+w_2\cdot r^i_{\text{delay}}+w_3\cdot r^i_{\text{speed}}+w_4\cdot r^i_{\text{rear}}+w_5\cdot r^i_{\text{lateral}},
\end{equation} 
where $w_1,\ldots,w_5\in\mathbb{R}_{\geq0}$ are the weighting factors corresponding to the following costs.
\paragraph{Fuel Efficiency}
We use the L$^2$-norm of the control input at each time step $k$ as a penalty to reduce the control effort, which decreases fuel consumption.  
 \begin{equation}\label{eq:delay}
r^i_{\text{fuel}}= - \dfrac{\| u_k^i\|^2}{\max\{\|u^{\max}\|,\|u^{\min}\|\}}.
\end{equation} 
\paragraph{Delay}
To improve the travel time, we define the time delay as a difference between the travel time that CAV $i$ takes from the entry of the control zone until the exit of the control zone and the time it would take for CAV $i$ from the entry until the exit of the control zone with the initial speed. 
Considering CAV $i$ at time step $k$, the traveled distance measured from the entry of the control zone is $p_{k}^i$, and its entry speed is denoted by $v_0^i$, the normalized penalty corresponding to delay is
 \begin{equation}\label{eq:delay}
r^i_{\text{delay}}= - \dfrac{\left(k\Delta t-\dfrac{p_{k}^i}{v_0^i}\right)}{\dfrac{p_{k}^i}{v_0^i}},
\end{equation}
where $\Delta t\in\mathbb{R}_{>0}$ is the time step.
\paragraph{Speed Limits Violation}
For each CAV $i$, at each time step $k$, the speed is bounded by 
\begin{equation}\label{vconstraint}
    0\leq v^{\min}\leq v_{k}^i\leq v^{\max},
\end{equation}
where $v^{\min},v^{\max}$ are the minimum and maximum speed limit, respectively. To ensure the speed constraint does not become active, we have 
 \begin{equation}\label{eq:speedLimit}
    r^i_{\text{speed}} = \begin{cases}
             -1,& \text{if speed violates the constraint},
            \\
            0,& \text{otherwise}.
          \end{cases}
\end{equation}

\paragraph{Rear-end Safety} To ensure the absence of rear-end collision between CAV $i$ and a preceding CAV $j$ at time step $k$, we impose the following rear-end safety constraint
\begin{equation}\label{RearEndCons}
 p_{k}^j-p_{k}^i\geq d_{\text{safe}}, 
\end{equation}
where $d_{\text{safe}}\in\mathbb{R}_{>0}$ is a safe constant distance.
The associated penalty for violating rear-end safety at each time step $k$ is
 \begin{equation}\label{eq:speedLimit}
    r^i_{\text{rear}} = \begin{cases}
             -100,& \text{if}~~ p_{k}^j-p_{k}^i< d_{\text{safe}},
            \\
            0,&\text{if}~~  p_{k}^j-p_{k}^i\geq d_{\text{safe}}.
          \end{cases}
\end{equation} 
\paragraph{Lateral Safety}
To guarantee lateral safety as CAVs cross the merging zone, we limit the merging zone occupancy to only one CAV at a time for CAVs with lateral collision potentials. For instance, in Fig. \ref{fig:1}, CAV \#$7$ or CAV \#$6$ is not allowed to be inside the merging zone at the same time. On the other hand, CAV \#$6$, CAV \#$9$, and CAV \#$5$ can be inside the merging zone at the same time (recall that right/left turns are not allowed).
Considering that CAVs $i$ and $j$ might have a lateral collision inside the merging zone, we construct the penalty according to the following two cases:

Case 1: CAV $i$ is outside the merging zone at time step $k$, and by taking an action $u^i_k$, it enters the merging zone while  CAV $j$, which previously entered the merging zone, is either still inside the merging zone, or it enters the merging zone at the same time as CAV $i$. In this case, CAV $i$ receives the  penalty $r^i_{\text{lateral}}=-100$.

Case 2: CAV $i$ is outside the merging zone at time step $k$, and by taking an action $u^i_k$, it enters the merging zone while, at the same time step $k$, CAV $j$ exits the merging zone. In this case, CAV $i$ receives no penalty.

\begin{figure*}[t]
  \subfloat[][]{\includegraphics[width=0.2\linewidth]{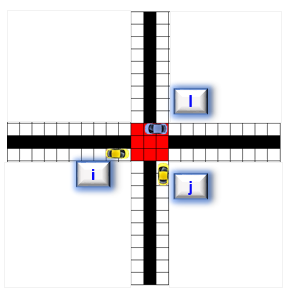}\label{a}}\quad
  \subfloat[][]{\includegraphics[width=0.2\linewidth]{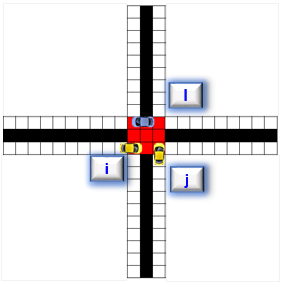}\label{b}}\quad\quad~~
 \subfloat[][]{\includegraphics[width=0.2\linewidth]{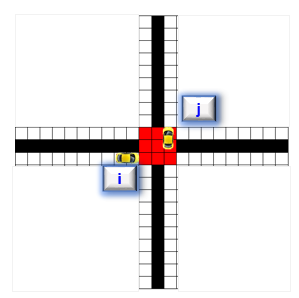}\label{c}}\quad
  \subfloat[][]{\includegraphics[width=0.2\linewidth]{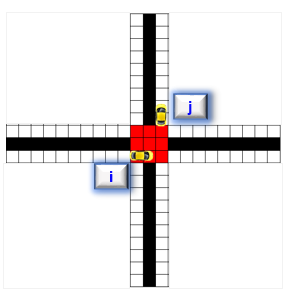}\label{d}}\quad
\caption{Case 1 in the lateral safety : \protect\subref{a} CAV $i$ and CAV $j$ are outside and CAV $l$ is inside  the merging zone at time step $k$; \protect\subref{b} CAV $i$, CAV $j$, and CAV $l$ are all inside at time step $k+1$; Case 2 in the lateral safety \protect\subref{c} CAV $i$ is outside and CAV $j$ is inside the merging zone at time step $k$; \protect\subref{d} CAV $i$ is inside and CAV $j$ is outside the merging zone at time step $k+1$.}\label{fig1:lateralCases} \end{figure*}

If by the time CAV $i$ enters the merging zone, there are more than one vehicles inside the merging zone that can cause lateral collision with CAV $i$, then CAV $i$ gets penalized for each of these CAVs separately, as follows. 
For the cases illustrated in Fig. \ref{fig1:lateralCases}\protect\subref{a} and \ref{fig1:lateralCases}\protect\subref{b} CAV $i$'s reward is $-100$ while CAVs $j$ and $l$ get $-200$ and $0$ reward, respectively; in Fig. \ref{fig1:lateralCases}\protect\subref{c} and \ref{fig1:lateralCases}\protect\subref{d} CAVs $i$ and $j$ both get reward $0$. To encourage CAVs to have a safe pass through the intersection, each CAV receives a terminal reward equal to $10n$ (recall that $n$ is the total number of CAVs) when it exits the control zone, if the episode did not have any crashes. Additionally, the episodes with a crash are terminated, and a new episode starts.  Algorithm  \ref{Alg:Q-learning} shows the pseudo-code of our decentralized coordination-aware framework for CAV $i$. 
\begin{algorithm}
 \caption{Hysteretic Q-learning algorithm for CAV $i$ }
 \hspace*{\algorithmicindent} \textbf{Input:} $\alpha,~\beta$  \\

 \begin{algorithmic}[1]
 \State {Initialize $Q^i(s_k^i,u_k^i)=0~~ \forall s_k^i\in\mathcal{S}^{i},u_k^i\in\mathcal{U}^{i} $}
\For {$t=$ 1 : total episodes}
\State {Reset entry time and speed randomly}
\For {$k=$ 1: maximum time step}
\If {CRASHED}
\State {$Break$}
\EndIf 
\State {$s_k^i \leftarrow$ Compute current state }
\State {$u_k^i \leftarrow$ Pick an epsilon-greedy action}
\State {$s_{k+1}^i \leftarrow$ Simulation }

\State {CRASHED $\leftarrow$ Did any CAVs violate the rear-end or lateral safety constraint?}
\State {$r_{k}^i \leftarrow$ Compute reward }
\State {$\delta \leftarrow r_k^i + \gamma \max\limits_{u_{k+1}^i} Q^i(s_{k+1}^i,u_{k+1}^i) -Q^i(s_k^i,u_k^i)$}

\If{$\delta\geq 0$}
 \State {$Q^i(s_k^i,u_k^i)\leftarrow Q^i(s_k^i,u_k^i) + \alpha \delta$}
\Else
 \State {$Q^i(s_k^i,u_k^i)\leftarrow Q^i(s_k^i,u_k^i) + \beta \delta$}
\EndIf
\EndFor
\EndFor
 \end{algorithmic} \label{Alg:Q-learning}
 \end{algorithm}

\subsection{FIFO queuing policy}
In this subsection, we provide a brief overview of a common approach in motion planning of CAVs at signal-free intersections called FIFO queuing policy. By imposing a FIFO queuing policy, each CAV must enter the merging zone in the same order that it entered the control zone. Let $\mathcal{N}(t)=\{1,\ldots,N(t)\}$ be the queue which designates the order in which CAVs enter the merging zone, and $t_i^0,~ t_i^m\in\mathbb{R}_{>0}$ denote the time when CAV $i\in\mathcal{N}(t)$ enters the control zone and merging zone, respectively.
The optimal entry time, which satisfies safety and speed constraint, can be found through the following recursive structure \cite{Malikopoulos2017}

\begin{equation}\label{tbar}
    t_i^{m^\ast} = \begin{cases}
            t_i^0+\dfrac{L}{v_0^i}  ,&\text{if}\;\: i=1,\\
            \\
           \max\{t_{i-1}^{m^\ast},~t_{j}^{m^\ast}+t_h,~t_i^c\}  ,&\text{if}\;\: i-1 \in\text{\textit{safe}},\\
            \\
            \max\{t_{i-1}^{m^\ast}+t_h,~t_i^c\}  ,&\text{if}\;\: i-1 \in\text{ \textit{lateral}},\\
            \\
            \max\{t_{i-1}^{m^\ast}+t_h,~t_i^c\}  ,&\text{if}\;\: i-1 \in\text{ \textit{rear-end}},\\
            \\
          \end{cases}
\end{equation} 
 where based on the path of CAV $i-1$ and CAV $i$, CAV $i-1$ belongs to one of the following subsets: (1) \textit{safe}, if there is no potential for collision with CAV $i$. (2) \textit{lateral}, if there is a potential for lateral collision with CAV $i$. (3) \textit{rear-end}, if CAV $i-1$ is the CAV immediately positioned in front of CAV $i$. The earliest feasible time that CAV $i$ can reach the merging zone is denoted by $t_i^c$, and $t_h$ is the safe time-headway to ensure safety at the entrance of the merging zone. In the above equation, if $i=1$, CAV $i$ cruises with the constant speed that it entered the control zone. Index $j$ in $t_j^m$ represents CAV $j$ which is physically located in front of CAV $i$.  
 \subsection{Combined Hysteretic Q-learning with FIFO Framework}
 In our hysteretic Q-learning framework (Algorithm \ref{Alg:Q-learning}), we aimed at achieving the lateral safety through our state and reward architecture. Due to the fact that after each crash the simulation episode is terminated, an increasing number of CAVs may require a greater number of simulation episodes, which might become less applicable in the real systems. 
 In this subsection, we introduce an enhanced framework which is a combination of FIFO and hysteretic Q-learning. In this framework, CAVs first find the optimal arrival time at merging zone recursively through a FIFO queuing policy at the start of each simulation episode. Since the lateral safety, and time-delay minimization are considered in the FIFO queuing policy \cite{Malikopoulos2017}, we need to modify the state and reward function in our Q-learning framework. The revised state of CAV $i$ at time step $k$ is $s^i_k\coloneqq\left\langle p^i_k,~v^i_k,~\mathcal{X}^{\text{i,rear}}_k,\Delta t_i^{m^\ast}\right\rangle$, where $p^i_k$ and $v^i_k$ are its position and speed, respectively; and $\mathcal{X}^{\text{i,rear}}_k\coloneqq\left\langle p^j_k,~v^j_k\right\rangle$ consists of the position and speed of CAV $j$, which is immediately ahead of CAV $i$. The difference between the optimal arrival time at merging zone, and arrival time at the control zone is captured in the last element $\Delta t_i^{m^\ast} = t_i^{m^\ast}-t_i^0$, which takes value from a bounded set defined by speed limits of the roads. In the revised reward function, the weights regarding the lateral collision and delay terms are set to zero
\begin{equation}\label{eq:modifedReward}
r_k^i=w_1^{\prime}\cdot r^i_{\text{fuel}}+w_2^{\prime}\cdot r^i_{\text{speed}}+w_3^{\prime}\cdot r^i_{\text{rear}}+w_4^{\prime}\cdot r^i_{\text{FIFO}},
\end{equation} 
where $w_1^{\prime},\ldots,w_4^{\prime}\in\mathbb{R}$ are new weighting factors.
To encourage CAV $i$ to reach the merging zone at the planned arrival time $t_i^{m^{\ast}}$, we define $r^i_{\text{FIFO}}$ to be the negative of the normalized squared error of arrival time at the merging zone computed as the difference of the estimated arrival time (EAT) at the merging zone and the optimal arrival time. At time step $k$, the EAT of CAV $i$ is approximated by assuming that CAV $i$ cruises with a constant speed $v^i_{k}$ for the rest of the path until the merging zone.
\begin{equation}\label{eq:fifo2}
r^i_{\text{FIFO}}= - ( \text{EAT} - t_i^{m^\ast})^2.
\end{equation}
As CAV $i$ enters the merging zone, the crossing time $t_i^m$ is compared to the optimal arrival time $t_i^{m^\ast}$, and CAV $i$ receives the last FIFO reward as follows 
\begin{equation}\label{eq:fifo3}
r^i_{\text{FIFO}}= - 10\times(t_i^m - t_i^{m^\ast})^2.
\end{equation}
After entering the merging zone, the corresponding weight ($w_4^{\prime}$) for CAV $i$ is set to zero.

\section{Simulation Setup}

\subsection{Discretization}
In our decentralized hysteretic Q-learning framework, at each time step each CAV needs to store the updated Q-value corresponding to the pair of current state and selected action. The discretization level not only directly affects the size of the Q-table, which each CAV stores, but it also influences the performance of our approach. Hence, determining proper discretization levels is a trade-off between the Q-table size and performance of the algorithm. Moreover, selecting improperly large or small values for time discretization results in poor performance, or even oscillating behavior. For instance, selecting a very small time step compared to the state discretization level leads to a situation in which a CAV takes action, but its state does not change. On the other hand, by selecting a very large time step, the system's safety might be jeopardized. The discretization level for position, speed, control input, and time are denoted by $\Delta p$, $\Delta v$, $\Delta u$, and $\Delta t$, respectively, as shown in Table \ref{tbl:Discretization level}.

\begin{table}[htbp]
\caption{Discretization level}
\vspace{0.5em}
\centering
\begin{tabular}{c|c|c|c|c} \label{tbl:Discretization level}
    Parameter & $\Delta p$& $\Delta v$ & $\Delta u$ & $\Delta t$
    \\
      (unit) & (m) & (m/s) & (m/s$^2$)& (s)
    \\
    \toprule
        Value &$2$&  $5$  &   $1$  & $0.5$\\
\end{tabular}
\end{table}

The other important parameters in our framework include those involved in the action selection determining the exploration and exploitation rate and the learning rates in the hysteretic update which are the step size of each Q-table update. These parameters are listed in Table \ref{tbl:simulationparameters}.

\begin{table}[htbp]
\caption{Selected parameters}
\vspace{0.5em}
\centering
\begin{tabular}{c|c|c|c|c} \label{tbl:simulationparameters}
    Parameter & $\epsilon_i$& $\epsilon_f$ & $\alpha$ & $\beta$
    \\
    \toprule
        Value &$0.6$&  $0.01$  &   $0.4$  & $0.05$\\
\end{tabular}
\end{table}

At the start of each episode, the initial conditions of CAVs are reset. In order to do that, the initial speed is drawn randomly from a uniform feasible speed distribution, and the arrival time of CAV $i$ is computed as 
\begin{equation}\label{eq:fifo2}
t_i^0 = \sum_{a=0}^{i} Y_a,
\end{equation}
where $t_i^0$ is the sum of $i$ independent and identically distributed random variable $Y_a$ drawn from an exponential distribution with mean $2$ s. 

\section{Simulation Results}
To evaluate the effectiveness of our proposed framework, we investigate the coordination of CAVs at a signal-free intersection in two scenarios (see Fig. \ref{fig:1}). We use the following parameters for the simulation: $d_{\text{safe}} = 4$ m, $v_{\min}=5$ m/s, $v_{\max}=15$ m/s, $u_{\max}=3$ m/s$^2$, $u_{\min}=-3$ m/s$^2$. 

\begin{figure}
\centering
\includegraphics[width=0.95\linewidth]{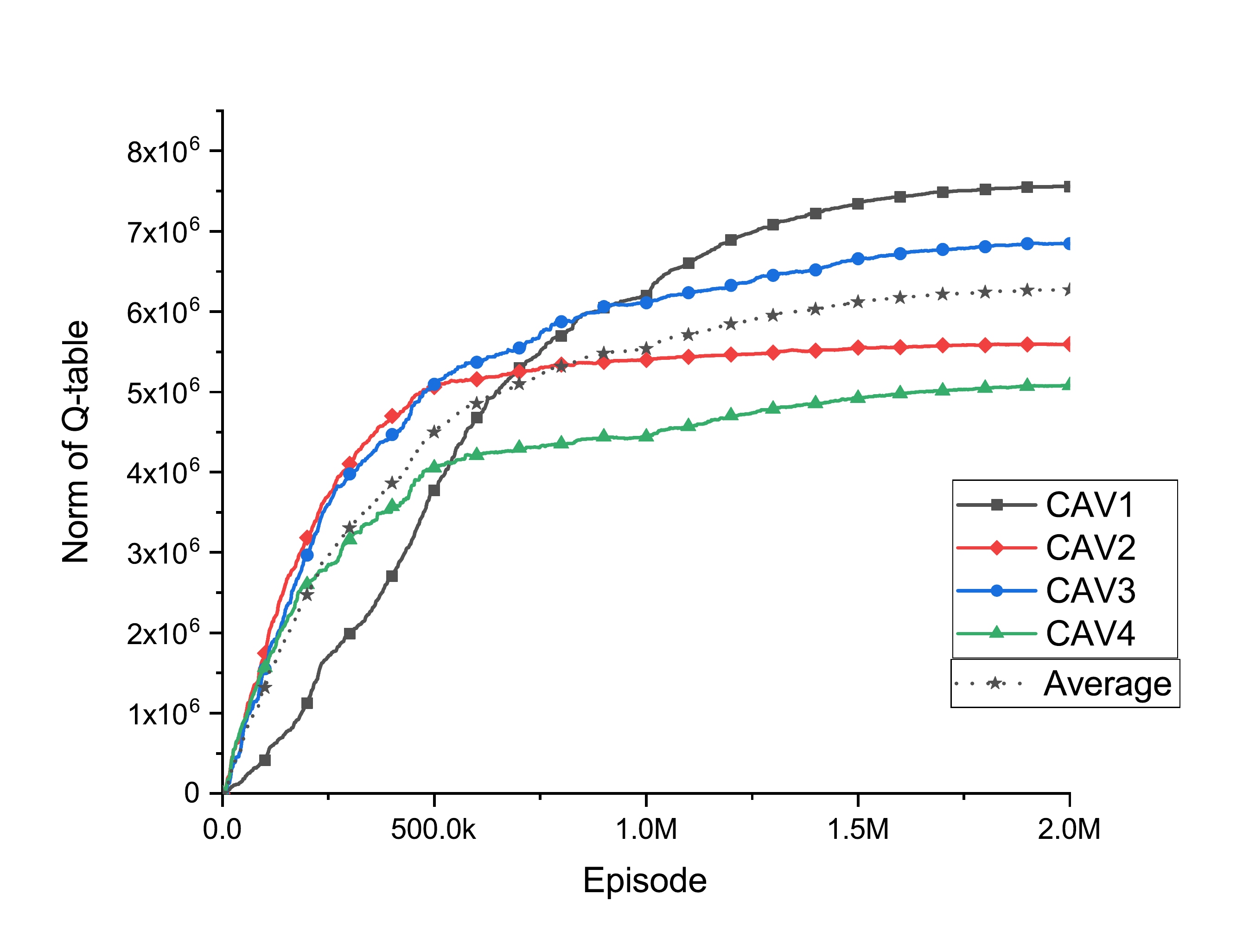}
\caption{Norm of Q-table for CAVs in Scenario 1. }
\label{fig:Q-tableNorm}%
\end{figure}

\textit{Scenario 1:} For our first scenario, we consider coordination of four CAVs using the hysteretic Q-learning framework in an intersection where the length of each road connected to the intersection is $L=32$ m, the length of the merging zones is $D=18$ m, and total episodes of simulation are set to $\num{2000000}$. CAV \#$1$, \#$2$, \#$3$, \#$4$, enter the control zone from southbound (SB), eastbound (EB), northbound (NB), and westbound (WB), respectively. Figure  \ref{fig:Q-tableNorm} shows the Cartesian norm of Q-table for each CAV and the average norm, respectively; which are computed after each $100$ episodes. As it can be seen in Fig. \ref{fig:Q-tableNorm}, the norm of Q-table for all four CAVs reach to stable values and does not change that much at the final episodes. The difference in the converged value is due to the fact that during the earlier episodes of training, CAV \#$3$ and \#$4$ are more likely to cause accidents and get penalized compared to CAV \#$1$ and CAV \#$2$, since CAV \#$3$ and \#$4$ enter the control zone later.

\begin{figure}
\centering
\includegraphics[width=0.95\linewidth]{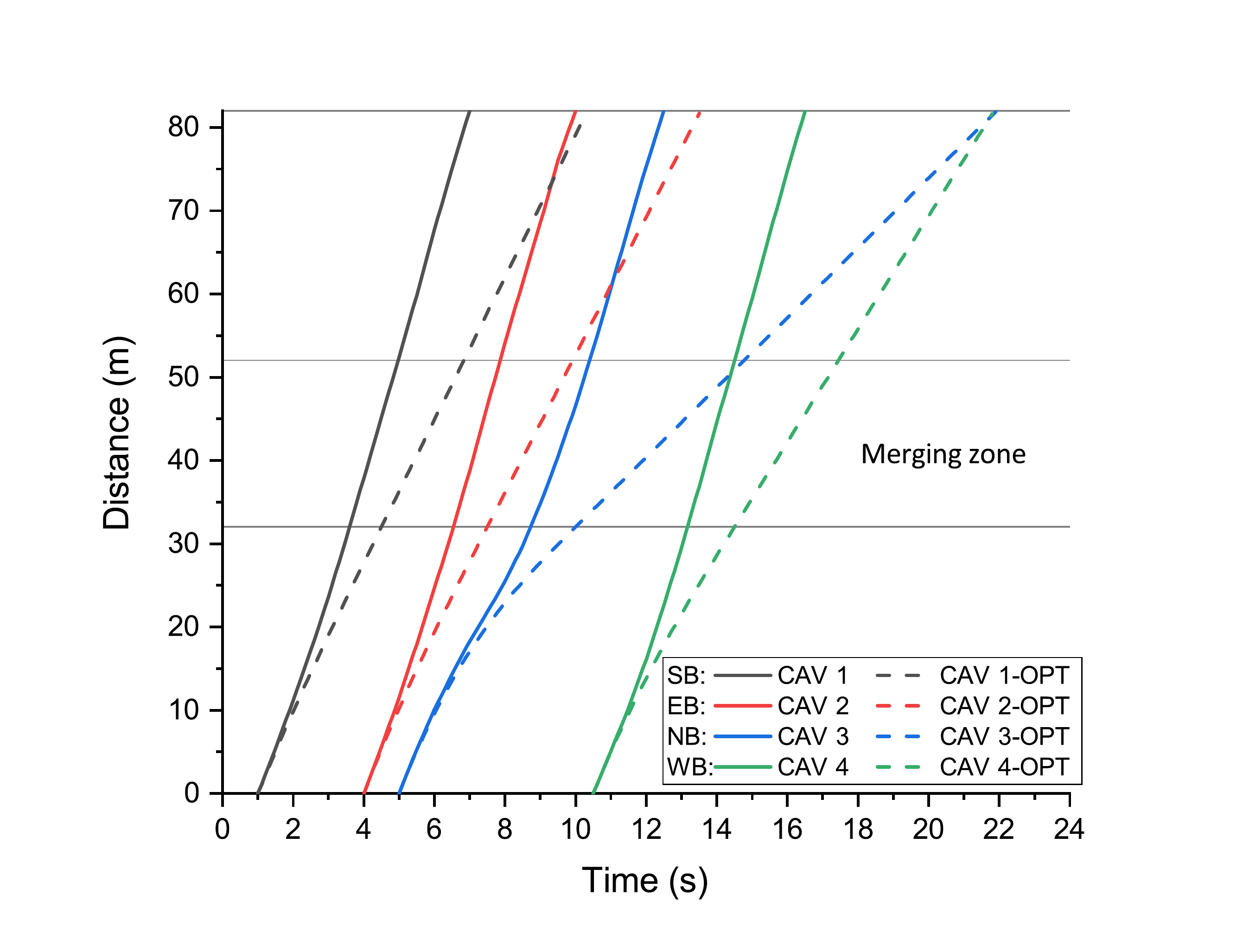}
\caption{Position trajectories of four trained CAVs and four CAVs found from the optimal control.}
\label{fig:PosTraj4}%
\end{figure}

After the training phase, we test the policies for $\num{1000}$ randomly generated simulation. The same initial conditions for four CAVs are used to simulate the optimal control framework. The optimal control framework consists of throughput maximization and energy minimization problems. In the throughput maximization problem, each CAV computes its arrival time at the merging zone based on a FIFO queuing policy. By restricting CAVs to have a constant speed after entering the merging zone, each CAV derives its energy-optimal control input from the control zone's entry until it reaches the merging zone considering speed and control constraints. Details of this approach can be found in \cite{Malikopoulos2017}. The position trajectories of four CAVs for a randomly selected simulation is shown in Fig. \ref{fig:PosTraj4} (solid lines). The major difference in the trajectories is because in the optimal control framework (dashed lines) the arrival time at merging zones for each CAV is found first, and then for each CAV, the optimal control problem is formulated from the arrival time at the control zone to the arrival time at the merging zone. On the other hand, our hysteretic Q-learning approach determines the policy with respect to the designed reward. Although trajectories resulted from our approach happen to respect the FIFO queuing policy (i.e., CAVs enter the merging zone in the same order they entered the control zone) without being enforced to, they appear to be more aggressive in minimizing the travel time. One can explore tuning $w_1$ and $w_2$ in \eqref{eq:rewardtotal} to find the trade-off between minimizing fuel consumption or delay.

\begin{figure}
\centering
\includegraphics[width=0.95\linewidth]{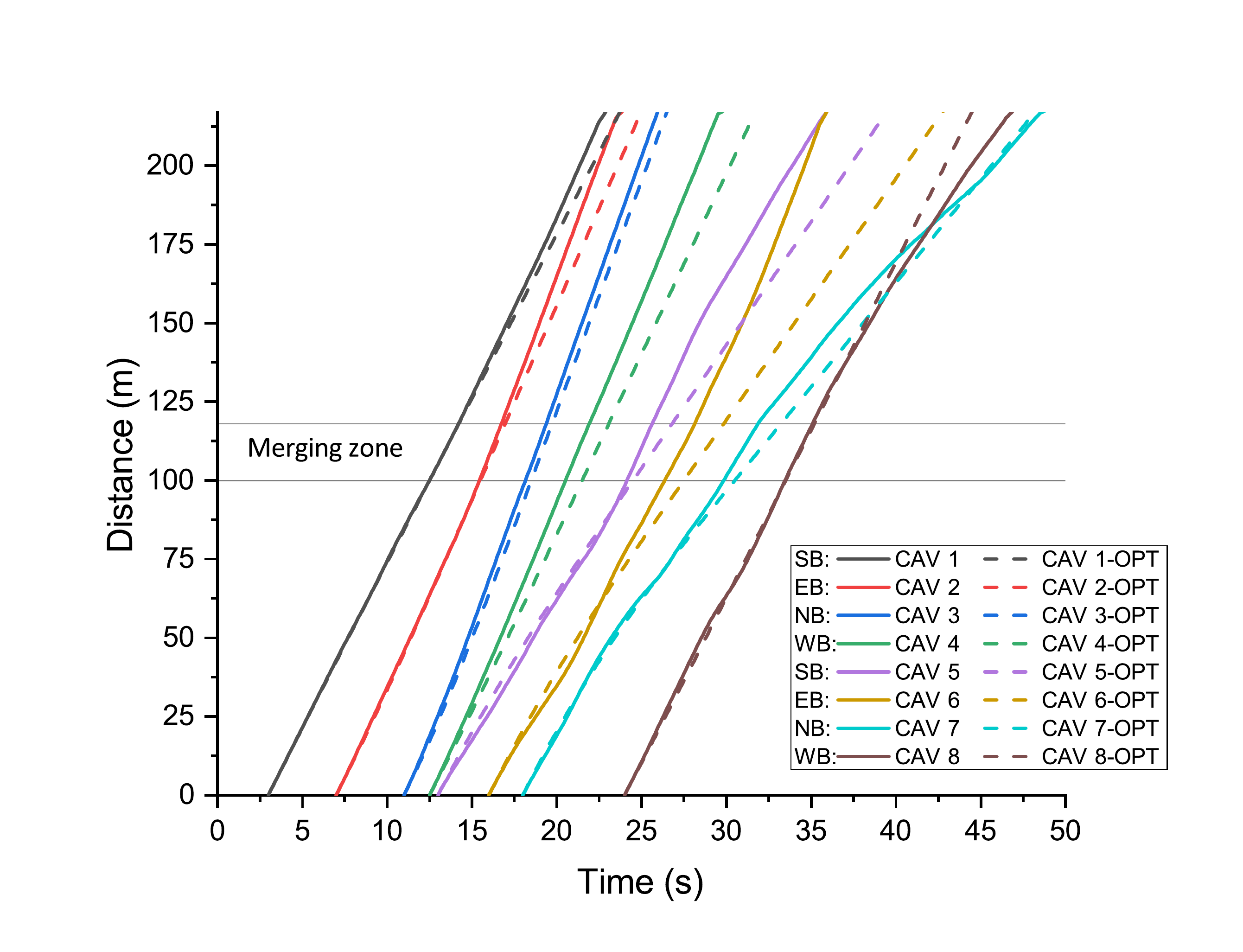}
\caption{Position trajectories of eight trained CAVs and eight CAVs found from the optimal control. }
\label{fig:PosTraj8}%
\end{figure}

\textit{Scenario 2:} For the second scenario, we consider coordination of eight CAVs using the combined hysteretic Q-learning with FIFO framework in an intersection which each road connecting to the intersection to be $L=100$ m, the length of the merging zones to be $D=18$ m, and total episodes of simulation are set to $\num{400000}$. In this scenario, we employ the state and reward architecture based on FIFO queuing policy. This extension allows us to reduce the state space size significantly. Namely, we are able to increase the control zone length in order to have enough space for CAVs to coordinate with each other. The position trajectories of eight CAVs for a randomly selected simulation (solid lines) along with the corresponding trajectories computed from the optimal control (dashed lines) are shown in Fig. \ref{fig:PosTraj8}. We note that CAVs following our combined hysteretic Q-learning with FIFO framework arrive at the merging zone at the planned arrival time with very small deviation. The trajectories for CAVs in our combined hysteretic Q-learning with FIFO framework do not deviate very much from the energy-optimal trajectories found from the optimal control techniques. Our proposed RL-based approach requires more time in the training phase compared to the classical control techniques, but after that Q-table is converged, it can be implemented in real-time. 
The videos from our simulation analysis can be found at the supplemental site, \url{https://sites.google.com/view/ud-ids-lab/HQLC}

\section{Concluding Remarks and Discussion} 
In this paper, we proposed a learning-based decentralized coordination framework for CAVs at a signal-free intersection to minimize travel delay and improve fuel consumption while ensuring rear-end and lateral safety. We embedded a coordination mechanism into our decentralized learning framework by using hysteretic Q-learning to update the Q-table of each CAV. We also integrated FIFO queuing policy in our framework to improve the performance of our system. Finally, we showed the effectiveness of our proposed approach through simulation and comparison with the classical optimal control method based on Pontryagin's minimum principle.

Ongoing research considers the presence of noise in the framework originated from the vehicle-level control and also investigates the effects of errors and delays in the communication. Our framework can be further extended to include lane-changing maneuvers and left/right turns by considering a different queuing policy instead of FIFO, such as upper-level motion planning proposed in \cite{Malikopoulos2020}. Coordination for mixed-traffic scenarios and the interaction of human-driven vehicles and CAVs, is another potential direction for future research. Future studies should also investigate approaches for transferring the policy to real-world scenarios. We have explored zero-shot transfer of an autonomous driving policy inside a roundabout directly from simulation to University of Delaware's scaled smart city under Gaussian noise \cite{jang2019simulation} and multi-agent adversarial noise \cite{chalaki2020ICCA}.

\bibliographystyle{IEEEtran.bst} 
\bibliography{reference/IDS_Publications_10302020.bib, reference/ref.bib}

\begin{thebibliography}{10}
\providecommand{\url}[1]{#1}
\csname url@rmstyle\endcsname
\providecommand{\newblock}{\relax}
\providecommand{\bibinfo}[2]{#2}
\providecommand\BIBentrySTDinterwordspacing{\spaceskip=0pt\relax}
\providecommand\BIBentryALTinterwordstretchfactor{4}
\providecommand\BIBentryALTinterwordspacing{\spaceskip=\fontdimen2\font plus
\BIBentryALTinterwordstretchfactor\fontdimen3\font minus
  \fontdimen4\font\relax}
\providecommand\BIBforeignlanguage[2]{{%
\expandafter\ifx\csname l@#1\endcsname\relax
\typeout{** WARNING: IEEEtran.bst: No hyphenation pattern has been}%
\typeout{** loaded for the language `#1'. Using the pattern for}%
\typeout{** the default language instead.}%
\else
\language=\csname l@#1\endcsname
\fi
#2}}

\bibitem{Schrank2019}
B.~Schrank, B.~Eisele, and T.~Lomax, ``{2019 Urban Mobility Scorecard},'' Texas
  A\& M Transportation Institute, Tech. Rep., 2019.

\bibitem{FHWA1}
{Federal Highway Administration}, ``{Department of Transportation},''
  {\url{www.safety.fhwa.dot.gov/intersection}}, online.

\bibitem{Klein2016a}
I.~Klein and E.~Ben-Elia, ``{Emergence of cooperation in congested road
  networks using ICT and future and emerging technologies: A game-based
  review},'' \emph{Transportation Research Part C: Emerging Technologies},
  vol.~72, pp. 10--28, 2016.

\bibitem{Melo2017a}
S.~Melo, J.~Macedo, and P.~Baptista, ``{Guiding cities to pursue a smart
  mobility paradigm: An example from vehicle routing guidance and its traffic
  and operational effects},'' \emph{Research in Transportation Economics},
  vol.~65, pp. 24--33, 2017.

\bibitem{ersal2020connected}
T.~Ersal, I.~Kolmanovsky, N.~Masoud, N.~Ozay, J.~Scruggs, R.~Vasudevan, and
  G.~Orosz, ``Connected and automated road vehicles: state of the art and
  future challenges,'' \emph{Vehicle system dynamics}, vol.~58, no.~5, pp.
  672--704, 2020.

\bibitem{Wadud2016}
Z.~Wadud, D.~MacKenzie, and P.~Leiby, ``Help or hindrance? the travel, energy
  and carbon impacts of highly automated vehicles,'' \emph{Transportation
  Research Part A: Policy and Practice}, vol.~86, pp. 1--18, 2016.

\bibitem{Levine1966}
W.~Levine and M.~Athans, ``{On the optimal error regulation of a string of
  moving vehicles},'' \emph{IEEE Transactions on Automatic Control}, vol.~11,
  no.~3, pp. 355--361, 1966.

\bibitem{Athans1969}
M.~Athans, ``{A unified approach to the vehicle-merging problem},''
  \emph{Transportation Research}, vol.~3, no.~1, pp. 123--133, 1969.

\bibitem{Malikopoulos2017}
A.~A. Malikopoulos, C.~G. Cassandras, and Y.~Zhang, ``A decentralized
  energy-optimal control framework for connected automated vehicles at
  signal-free intersections,'' \emph{Automatica}, vol.~93, pp. 244--256, 2018.

\bibitem{hult2018optimal}
R.~Hult, M.~Zanon, S.~Gros, and P.~Falcone, ``Optimal coordination of automated
  vehicles at intersections: Theory and experiments,'' \emph{IEEE Transactions
  on Control Systems Technology}, vol.~27, no.~6, pp. 2510--2525, 2018.

\bibitem{bichiou2018developing}
Y.~Bichiou and H.~A. Rakha, ``Developing an optimal intersection control system
  for automated connected vehicles,'' \emph{IEEE Transactions on Intelligent
  Transportation Systems}, vol.~20, no.~5, pp. 1908--1916, 2018.

\bibitem{Lee2012a}
J.~Lee and B.~Park, ``Development and evaluation of a cooperative vehicle
  intersection control algorithm under the connected vehicles environment,''
  \emph{IEEE Transactions on Intelligent Transportation Systems}, vol.~13,
  no.~1, pp. 81--90, 2012.

\bibitem{chalaki2020TCST}
B.~Chalaki and A.~A. Malikopoulos, ``Optimal control of connected and automated
  vehicles at multiple adjacent intersections,'' \emph{arXiv preprint
  arXiv:2008.02379}, 2020.

\bibitem{Ntousakis:2016aa}
I.~A. Ntousakis, I.~K. Nikolos, and M.~Papageorgiou, ``Optimal vehicle
  trajectory planning in the context of cooperative merging on highways,''
  \emph{Transportation Research Part C: Emerging Technologies}, vol.~71, pp.
  464--488, 2016.

\bibitem{xiao2019decentralized}
W.~Xiao, C.~Belta, and C.~G. Cassandras, ``Decentralized merging control in
  traffic networks: A control barrier function approach,'' in \emph{Proceedings
  of the 10th ACM/IEEE International Conference on Cyber-Physical Systems},
  2019, pp. 270--279.

\bibitem{Rios-Torres2017b}
J.~Rios-Torres and A.~A. Malikopoulos, ``{Automated and Cooperative Vehicle
  Merging at Highway On-Ramps},'' \emph{IEEE Transactions on Intelligent
  Transportation Systems}, vol.~18, no.~4, pp. 780--789, 2017.

\bibitem{Malikopoulos2018c}
A.~A. {Malikopoulos}, S.~{Hong}, B.~B. {Park}, J.~{Lee}, and S.~{Ryu},
  ``Optimal control for speed harmonization of automated vehicles,'' \emph{IEEE
  Transactions on Intelligent Transportation Systems}, vol.~20, no.~7, pp.
  2405--2417, 2019.

\bibitem{guanetti2018control}
J.~Guanetti, Y.~Kim, and F.~Borrelli, ``Control of connected and automated
  vehicles: State of the art and future challenges,'' \emph{Annual Reviews in
  Control}, vol.~45, pp. 18--40, 2018.

\bibitem{Rios-Torres2017}
J.~Rios-Torres and A.~A. Malikopoulos, ``A survey on the coordination of
  connected and automated vehicles at intersections and merging at highway
  on-ramps,'' \emph{IEEE Transactions on Intelligent Transportation Systems},
  vol.~18, no.~5, pp. 1066--1077, 2016.

\bibitem{sutton2018reinforcement}
R.~S. Sutton and A.~G. Barto, \emph{Reinforcement learning: An
  introduction}.\hskip 1em plus 0.5em minus 0.4em\relax MIT press, 2018.

\bibitem{kiumarsi2017optimal}
B.~Kiumarsi, K.~G. Vamvoudakis, H.~Modares, and F.~L. Lewis, ``Optimal and
  autonomous control using reinforcement learning: A survey,'' \emph{IEEE
  transactions on neural networks and learning systems}, vol.~29, no.~6, pp.
  2042--2062, 2017.

\bibitem{watkins1989learning}
C.~J. C.~H. Watkins, ``Learning from delayed rewards,'' 1989.

\bibitem{el2010agent}
S.~El-Tantawy and B.~Abdulhai, ``An agent-based learning towards decentralized
  and coordinated traffic signal control,'' in \emph{13th International IEEE
  Conference on Intelligent Transportation Systems}.\hskip 1em plus 0.5em minus
  0.4em\relax IEEE, 2010, pp. 665--670.

\bibitem{ji2009optimal}
X.~Ji and Z.~He, ``An optimal control method for expressways entering ramps
  metering based on q-learning,'' in \emph{2009 Second International Conference
  on Intelligent Computation Technology and Automation}, vol.~1.\hskip 1em plus
  0.5em minus 0.4em\relax IEEE, 2009, pp. 739--741.

\bibitem{jacob2005integrated}
C.~Jacob and B.~Abdulhai, ``Integrated traffic corridor control using machine
  learning,'' in \emph{2005 IEEE International Conference on Systems, Man and
  Cybernetics}, vol.~4.\hskip 1em plus 0.5em minus 0.4em\relax IEEE, 2005, pp.
  3460--3465.

\bibitem{davarynejad2011motorway}
M.~Davarynejad, A.~Hegyi, J.~Vrancken, and J.~van~den Berg, ``Motorway
  ramp-metering control with queuing consideration using q-learning,'' in
  \emph{2011 14th International IEEE Conference on Intelligent Transportation
  Systems (ITSC)}.\hskip 1em plus 0.5em minus 0.4em\relax IEEE, 2011, pp.
  1652--1658.

\bibitem{ivanjko2015ramp}
E.~Ivanjko, D.~K. Ne{\v{c}}oska, M.~Greguri{\'c}, M.~Vuji{\'c},
  G.~Jurkovi{\'c}, and S.~Mand{\v{z}}uka, ``Ramp metering control based on the
  q-learning algorithm,'' \emph{Cybernetics and Information Technologies},
  vol.~15, no.~5, pp. 88--97, 2015.

\bibitem{wang2019q}
L.~Wang, F.~Ye, Y.~Wang, J.~Guo, I.~Papamichail, M.~Papageorgiou, S.~Hu, and
  L.~Zhang, ``A q-learning foresighted approach to ego-efficient lane changes
  of connected and automated vehicles on freeways,'' in \emph{2019 IEEE
  Intelligent Transportation Systems Conference (ITSC)}.\hskip 1em plus 0.5em
  minus 0.4em\relax IEEE, 2019, pp. 1385--1392.

\bibitem{ngai2011multiple}
D.~C.~K. Ngai and N.~H.~C. Yung, ``A multiple-goal reinforcement learning
  method for complex vehicle overtaking maneuvers,'' \emph{IEEE Transactions on
  Intelligent Transportation Systems}, vol.~12, no.~2, pp. 509--522, 2011.

\bibitem{wu2019dcl}
Y.~Wu, H.~Chen, and F.~Zhu, ``Dcl-aim: Decentralized coordination learning of
  autonomous intersection management for connected and automated vehicles,''
  \emph{Transportation Research Part C: Emerging Technologies}, vol. 103, pp.
  246--260, 2019.

\bibitem{isele2017go}
D.~Isele and A.~Cosgun, ``To go or not to go: a case for q-learning at
  unsignalized intersections,'' 2017.

\bibitem{seliman2020automated}
S.~M. Seliman, A.~W. Sadek, and Q.~He, ``Automated vehicle control at freeway
  lane-drops: a deep reinforcement learning approach,'' \emph{Journal of Big
  Data Analytics in Transportation}, pp. 1--20, 2020.

\bibitem{tram2018learning}
T.~Tram, A.~Jansson, R.~Gr{\"o}nberg, M.~Ali, and J.~Sj{\"o}berg, ``Learning
  negotiating behavior between cars in intersections using deep q-learning,''
  in \emph{2018 21st International Conference on Intelligent Transportation
  Systems (ITSC)}.\hskip 1em plus 0.5em minus 0.4em\relax IEEE, 2018, pp.
  3169--3174.

\bibitem{meng2020adaptive}
L.~Meng-Lin, C.~Shao-Fei, and C.~Jing, ``Adaptive learning: A new decentralized
  reinforcement learning approach for cooperative multiagent systems,''
  \emph{IEEE Access}, 2020.

\bibitem{matignon2007hysteretic}
L.~Matignon, G.~J. Laurent, and N.~Le~Fort-Piat, ``Hysteretic q-learning: an
  algorithm for decentralized reinforcement learning in cooperative multi-agent
  teams,'' in \emph{2007 IEEE/RSJ International Conference on Intelligent
  Robots and Systems}.\hskip 1em plus 0.5em minus 0.4em\relax IEEE, 2007, pp.
  64--69.

\bibitem{Malikopoulos2020}
A.~A. Malikopoulos, L.~E. Beaver, and I.~V. Chremos, ``Optimal path planning
  and coordination for connected and automated vehicles,''
  \emph{arXiv:2003.12183}, 2020.

\bibitem{jang2019simulation}
K.~Jang, E.~Vinitsky, B.~Chalaki, B.~Remer, L.~Beaver, A.~A. Malikopoulos, and
  A.~Bayen, ``Simulation to scaled city: zero-shot policy transfer for traffic
  control via autonomous vehicles,'' in \emph{Proceedings of the 10th ACM/IEEE
  International Conference on Cyber-Physical Systems}, 2019, pp. 291--300.

\bibitem{chalaki2020ICCA}
B.~Chalaki, L.~E. Beaver, B.~Remer, K.~Jang, E.~Vinitsky, A.~Bayen, and A.~A.
  Malikopoulos, ``Zero-shot autonomous vehicle policy transfer: From simulation
  to real-world via adversarial learning,'' in \emph{IEEE 16th International
  Conference on Control \& Automation (ICCA)}, 2020, pp. 35--40.

\end{thebibliography}

\end{document}